\documentclass{amsart}
\usepackage{graphicx}
\usepackage{amsfonts,amssymb,amstext}
\usepackage{amsmath}
\usepackage{amstext}
\usepackage{amsthm}
\usepackage{amssymb,amsmath}
\usepackage{color}
\usepackage{graphics}

\usepackage[russian]{babel}
\usepackage{textcomp}
\usepackage{amsthm}
\usepackage{amsfonts}
\usepackage[mathcal]{eucal}
\usepackage{mathrsfs}

\numberwithin{equation}{section}
\newtheorem{theorem}{Теорема}[section]

\newtheorem{lemma}{Лемма}[section]
\newtheorem{remark}{Замечание}[section]

\makeatletter
\def\ps@firstpage{\ps@plain
\def\@oddfoot{\normalfont\normalsize \hfil\thepage\hfil
\global\topskip\normaltopskip}%
\let\@evenfoot\@oddfoot
\def\@oddhead{\@serieslogo\hss}%
\let\@evenhead\@oddhead % in case an article starts on a left-hand page
}

\def\ps@plain{\ps@empty
\def\@oddhead{\normalfont\scriptsize\hfil\MakeUppercase\shorttitle\hfil}
\def\@evenhead{\normalfont\scriptsize\hfil\MakeUppercase\shortauthors\hfil}
\def\@oddfoot{\normalfont\normalsize \hfil\thepage\hfil}%
\let\@evenfoot\@oddfoot
} \makeatother

\begin{document}
\pagestyle{plain}

\noindent{\sl Известия НАН Армении, Математика, том 52, н. 6, 2017, стр. 77-90}

\title[Сходимость разложений по собственным функциям ...]{Сходимость разложений по собственным
функциям и асимптотика спектральных данных задачи Штурма-Лиувилля}

\author{А. А. Пахлеванян}

\maketitle

\centerline{\sl Институт математики НАН Армении{\footnote{Исследование выполнено при финансовой поддержке
Государственного комитета по науке МОН РА в рамках научного проекта No. 16YR--1A017.}}}

\centerline{E-mail: {\it apahlevanyan@instmath.sci.am}}

\begin{abstract}
Доказывается равномерная сходимость разложения абсолютно непрерывной функции по собственным функциям
задачи Штурма-Лиувилля $-y'' + q \left( x \right) y = \mu y,$ $y \left(0\right)=0,$ $y\left( \pi \right)\cos
\beta + y'\left( \pi \right)\sin \beta  = 0,$ $\beta \in \left( 0, \pi \right)$ с суммируемым потенциалом
$q \in L_{\mathbb{R}}^1 \left[0, \pi \right].$ Этот результат используется для получения новых, более точных
асимптотических формул для собственных значений и нормировочных постоянных этой задачи.
\end{abstract}

\vskip2mm


\vskip2mm

\noindent{\bf Ключевые слова:} задача Штурма-Лиувилля; разложение по собственным функциям; асимптотика
собственных значений и нормировочных постоянных.

\baselineskip15pt

\section{Введение и формулировка результатов}
\setcounter{equation}{0}

Обозначим через $L\left(q, \alpha, \beta \right)$ следующую краевую задачу Штурма-Лиувилля
\begin{equation}\label{eq1.1}
- y'' + q\left( x \right)y = \mu y \equiv \lambda^2 y, \; x \in \left( 0, \pi \right), \; \mu \in \mathbb{C},
\end{equation}
\begin{equation}\label{eq1.2}
y\left( 0 \right)\cos \alpha  + y'\left( 0 \right)\sin \alpha  = 0, \; \alpha \in \left( 0, \pi \right],
\end{equation}
\begin{equation}\label{eq1.3}
y\left( \pi \right)\cos \beta  + y'\left( \pi \right)\sin \beta  = 0, \; \beta \in \left[ 0, \pi \right),
\end{equation}
где $q$ вещественная, суммируемая на $\left[0, \pi \right]$ функция (мы пишем $q \in L_{\mathbb{R}}^1\left[ 0, \pi \right]$). Через $L\left(q,\alpha, \beta \right)$ будем обозначать также самосопряженный оператор, порожденный задачей \eqref{eq1.1}--\eqref{eq1.3} в гильбертовом пространстве $L^2\left[0, \pi \right]$ (см. \cite{Naimark:1969,Levitan_Sargsyan:1970}). Хорошо известно, что спектр оператора $L\left(q, \alpha, \beta \right)$ чисто дискретен и состоит из простых, действительных собственных значений (см. \cite{Naimark:1969}---\cite{Marchenko:1977}), которые мы обозначаем через ${\mu}_n \left(q, \alpha, \beta \right),$ $n=0,1,2,\dots,$ подчеркивая зависимость от $q,$ $\alpha$ и $\beta.$ Предполагается, что собственные значения $\mu_n$ пронумерованы в порядке возрастания:
\begin{equation*}
\mu_0 \left(q, \alpha, \beta \right) < \mu_1 \left(q, \alpha, \beta \right) < \dots < \mu_n \left(q, \alpha, \beta \right) < \dots \; .
\end{equation*}
Обозначим через $\varphi \left(x, \mu, \alpha \right)$ и $\psi \left(x, \mu, \beta \right)$ решения уравнения \eqref{eq1.1}, удовлетворяющие начальным условиям
$$%\begin{equation}\label{eq1.4}
\varphi \left(0, \mu, \alpha \right)=\sin \alpha, \;\,\, \varphi' \left(0, \mu, \alpha \right)=-\cos \alpha,
$$%\end{equation}
$$%\begin{equation}\label{eq1.5}
\psi \left(\pi, \mu, \beta \right)=\sin \beta, \;\,\, \psi' \left(\pi, \mu, \beta \right)=-\cos \beta.
$$%\end{equation}
Хорошо известно (\cite{Naimark:1969,Levitan_Sargsyan:1970,Harutyunyan_Hovsepyan:2005}), что для фиксированного $x,$  $\varphi,$ $\varphi',$ $\psi,$ $\psi'$ являются целыми функциями от $\mu.$
Обозначим через $W_{\alpha, \beta}\left(x, \mu \right)$ вронскиан решений $\varphi \left(x, \mu, \alpha \right)$ и $\psi \left(x, \mu, \beta \right):$
\begin{equation}\label{eq1.6}
W_{\alpha, \beta}\left(x, \mu \right):=\varphi \left(x, \mu, \alpha \right) \psi' \left(x, \mu, \beta \right) - \varphi' \left(x, \mu, \alpha \right) \psi \left(x, \mu, \beta \right).
\end{equation}
Из формулы Лиувилля для вронскиана следует (см. например, \cite{Coddington_Levinson:1955}) что $W_{\alpha, \beta}\left(x, \mu \right)$ не зависит от $x.$ Легко видеть что функции ${\varphi_n}\left(x\right):=\varphi(x, \mu_n, \alpha)$ и ${\psi_n}\left(x\right):=\psi\left(x, \mu_n, \beta\right),$ $n=0,1,2,\dots,$ являются собственными функциями, соответствующими собственному значению $\mu_n.$ Так как собственные значения простые, то существуют числа $\beta_n=\beta_n\left(q, \alpha, \beta\right),$ $n=0,1,2,\dots,$ такие что
\begin{equation}\label{eq1.8}
\psi_n\left(x\right)=\beta_n \varphi_n\left(x\right), \beta_n \neq 0.
\end{equation}
Квадраты $L^2$-норм этих собственных функций:
\begin{equation}\label{eq1.9}
a_n=a_n \left(q, \alpha, \beta \right)=\int\limits_0^\pi {\varphi_n^2 \left(x\right)dx}, \quad b_n=b_n \left(q, \alpha, \beta \right)=\int\limits_0^\pi {\psi_n^2 \left(x\right)dx},
\end{equation}
называются нормировочными постоянными.

Одна из основных теорем спектральной теории дифференциальных операторов (см. \cite{Naimark:1969}) гласит:

\begin{theorem}(\cite[стр. 90]{Naimark:1969})\label{thm2.1}
Всякая функция из области определения самосопряженного дифференциального оператора разлагается в равномерно сходящийся обобщенный ряд Фурье по собственным функциям этого оператора.
\end{theorem}
Этот результат не может быть применен для функций, не принадлежащих области определения самосопряженного дифференциального оператора. В случае оператора Штурма-Лиувилля доказано (см. например, \cite{Yurko:2007}), что при условии $\sin \alpha \neq 0,$ $\sin \beta \neq 0$ (т.е. $\alpha, \beta \in \left(0, \pi \right)$) равномерная сходимость разложения имеет место для любой абсолютно непрерывной функции, а именно имеет место (см. также \cite{Coddington_Levinson:1955,Levitan_Sargsyan:1970,Levitan_Sargsyan:1988})

\begin{theorem}(\cite{Yurko:2007})\label{thm2.2}
Пусть $q \in L_{\mathbb{R}}^2 \left[0, \pi \right],$ $\alpha, \beta \in \left(0, \pi \right)$ и $f$ абсолютно непрерывная функция на $\left[0, \pi \right].$ Тогда
$$ \mathop {\lim}\limits_{N \to \infty}\mathop {\operatorname{max}}\limits_{x \in \left[0, \pi\right]}\left|f\left(x\right)-\sum\limits_{n = 0}^{N} c_n {\varphi_n \left(x \right)}\right|=0, \; c_n=\dfrac{1}{a_n}{\displaystyle\int_{0}^{\pi}f\left(t\right){\varphi_n \left(t \right)}dt},$$
где $\varphi_n \left( x \right) \equiv \varphi \left(x, \mu_n \left(q, \alpha, \beta \right), \alpha \right).$
\end{theorem}

Первой целью данной работы является доказательство того, что аналогичные результаты имеют место также для задач $L \left( q, \pi, \beta \right),$ $\beta \in \left(0,\pi \right)$ и $L \left( q, \alpha, 0 \right),$ $\alpha \in \left(0,\pi\right),$ причем при более общем условии $q \in L_{\mathbb{R}}^1 \left[0, \pi \right]:$

\begin{theorem}\label{thm2.3}
Пусть $q \in L_{\mathbb{R}}^1\left[ {0,\pi } \right],$ $\alpha=\pi,$ $\beta \in \left(0, \pi \right)$ и $f$ абсолютно непрерывная функция на $\left[0, \pi\right].$ Тогда для любого $a \in \left(0, \pi \right)$
\begin{equation}\label{eq2.1}
\mathop {\lim}\limits_{N \to \infty}\mathop {\operatorname{max}}\limits_{x \in \left[a, \pi\right]}\left|f\left(x\right)-\sum\limits_{n = 0}^{N} c_n {\varphi_n \left(x \right)}\right|=0, \; \,\, c_n=\dfrac{1}{a_n}{\displaystyle\int_{0}^{\pi}f\left(t\right){\varphi_n \left(t \right)}dt},
\end{equation}
где $\varphi_n \left( x \right) \equiv \varphi \left(x, \mu_n \left(q, \pi, \beta \right), \pi \right) \equiv \varphi \left(x, \mu_n, \pi \right).$
\end{theorem}

\begin{theorem}\label{thm2.4}
Пусть $q \in L_{\mathbb{R}}^1\left[ {0,\pi } \right],$ $\alpha \in \left(0, \pi \right),$ $\beta=0$ и $f$ абсолютно непрерывная функция на $\left[0, \pi\right].$ Тогда для любого $b \in \left(0, \pi \right)$
\begin{equation}\label{eq2.2}
\mathop {\lim}\limits_{N \to \infty}\mathop {\operatorname{max}}\limits_{x \in \left[0, b \right]}\left|f\left(x\right)-\sum\limits_{n = 0}^{N} c_n {\varphi_n \left(x \right)}\right|=0, \; c_n=\dfrac{1}{a_n}{\displaystyle\int_{0}^{\pi}f\left(t\right){\varphi_n \left(t \right)}dt},
\end{equation}
где $\varphi_n \left( x \right) \equiv \varphi \left(x, \mu_n \left(q, \alpha, 0 \right), \alpha \right).$
\end{theorem}

\begin{remark}\label{rem2.7}
Легко видеть, что улучшить результаты Теорем \ref{thm2.3} и \ref{thm2.4} (получить равномерную сходимость рядов в \eqref{eq2.1} и \eqref{eq2.2} на всем отрезке $\left[0, \pi \right]$ без дополнительных условий) невозможно. Действительно, для абсолютно непрерывной функции $f \equiv \dfrac{\pi}{2},$ $x \in \left[0, 2\pi \right],$ имеет место тождество (см. например, \cite[формула 37 на стр. 578]{Bronshtein_Semendyayev:1998})
\begin{equation*}
\dfrac{\pi}{2}=\displaystyle\sum\limits_{n = 0}^\infty \dfrac{\sin \left(n+\frac{1}{2}\right)x}{n+\frac{1}{2}}, \; 0 < x < 2\pi,
\end{equation*}
т.е. Теорема \ref{thm2.3} для нее верна, но если заменить $\mathop {\operatorname{max}}\limits_{x \in \left[a, \pi\right]}\left|\dots\right|$ на $\mathop {\operatorname{max}}\limits_{x \in \left[0, \pi\right]}\left|\dots\right|,$ то теорема перестает быть верной. С другой стороны, если в Теореме \ref{thm2.3} взять $f\left(0 \right)=0,$ тогда ряд в \eqref{eq2.1} сходится равномерно на всем отрезке $\left[0, \pi \right].$ Аналогичное утверждение верно для Теоремы \ref{thm2.4}, если взять $f\left(\pi \right)=0.$
\end{remark}
Второй целью нашей работы является получение асимптотической формулы для собственных значений задачи $L \left(q, \pi, \beta \right)$ при $q \in L_{\mathbb{R}}^1\left[0, \pi \right]$ и $\beta \in \left(0, \pi \right)$ (т.е. $\sin \beta \neq 0$). Прежде чем сформулировать результат, заметим, что в работе \cite{Harutyunyan:2008} Т. Н. Арутюнян ввел понятие функции $\delta_n \left(\alpha, \beta \right),$ которая определяется формулой
\begin{equation*}
\delta_n \left(\alpha, \beta \right):=\sqrt {\mu_n \left(0, \alpha, \beta \right)}-n=\lambda_n \left(0, \alpha, \beta \right)-\lambda_n \left(0, \dfrac{\pi}{2}, \dfrac{\pi}{2}\right), \; n \geq 2,
\end{equation*}
и доказал что $-1 \leq \delta_{n} \left(\alpha, \beta \right) \leq 1$ и $\delta_{n} \left(\alpha, \beta \right)$ является решением следующего трансцендентного уравнения:
\begin{multline}\label{eq1.10}
\delta_n \left(\alpha, \beta \right)=\dfrac{1}{\pi}\arccos \dfrac{\cos \alpha}{\sqrt {\left(n+{\delta_n} \left(\alpha, \beta \right)\right)^2 \sin^2 \alpha + \cos^2 \alpha}}- \\
-\dfrac{1}{\pi}\arccos \dfrac{\cos \beta}{\sqrt {\left(n+\delta_n \left(\alpha, \beta \right)\right)^2 \sin^2 \beta+\cos^2 \beta}}.
\end{multline}

\begin{theorem}\label{thm3.1}
Пусть $q \in L_{\mathbb{R}}^1 \left[0, \pi \right]$ и пусть $\lambda_n^2\left(q, \alpha, \beta\right)={\mu _n}\left(q, \alpha, \beta\right).$ Тогда
\begin{enumerate}
\item[(a)] имеет место асимптотическое соотношение $\left(n \to \infty \right)$
\begin{equation}\label{eq3.1}
\lambda_n\left(q, \alpha, \beta\right)=n+\delta_n\left(\alpha, \beta\right)+ \dfrac{[q]}{2\left(n+\delta_n\left(\alpha, \beta\right)\right)}+l_n\left(q, \alpha, \beta\right)+O\left(\dfrac{1}{n^2}\right),
\end{equation}
где $\left[q \right]=\dfrac{1}{\pi}\displaystyle\int_{0}^{\pi}q\left(t\right)dt,$
\begin{equation*}
l_n\left(q, \alpha, \beta\right)=\dfrac{1}{2 \pi \left(n+\delta_n\left(\alpha, \beta\right)\right)} \int_{0}^{\pi} q(x)\cos2\left(n+\delta_n\left(\alpha, \beta\right)\right)xdx, \; \alpha \in \left(0, \pi \right),
\end{equation*}
\begin{equation}\label{eq3.2}
l_n=l_n\left(q, \pi, \beta\right)=-\dfrac{1}{2 \pi \left(n+\delta_n\left(\pi, \beta\right)\right)} \int_{0}^{\pi} q(x)\cos2\left(n+\delta_n\left(\pi, \beta\right)\right)xdx.
\end{equation}
Оценка остатка $O\left(\dfrac{1}{n^2}\right)$ в \eqref{eq3.1} равномерна по всем $\alpha, \beta \in [0, \pi]$ и $q$ из ограниченных подмножеств $L_{\mathbb{R}}^1\left[ {0,\pi } \right]$ (мы будем писать $q \in {BL}^1_\mathbb{R}\left[0, \pi\right]$).
\item[(b)] Функция $l,$ определенная формулой
\begin{equation}\label{eq3.3}
l(x)=\displaystyle\sum_{n=2}^{\infty} l_n\left(q, \alpha, \beta\right) \sin \left(n+\delta_n\left(\alpha, \beta\right)\right)x,
\end{equation}
абсолютно непрерывна на произвольном отрезке $\left[a,b\right] \subset \left(0, 2 \pi\right),$ т.е. $l \in AC\left(0,2\pi\right).$
\end{enumerate}
\end{theorem}
В работе \cite{Harutyunyan:2016} утверждение (b) теоремы \ref{thm3.1} было доказано при условии $\alpha, \beta \in \left(0, \pi \right)$ и в случае $\alpha=\pi,$ $\beta=0.$ Мы докажем, что это утверждение верно также при $\alpha=\pi,$ $\beta \in \left(0, \pi \right).$ Этому посвящен раздел \ref{sec3}. Третьей целью нашей работы является получение асимптотических формул для нормировочных постоянных $a_n$ и $b_n$ (см. \eqref{eq1.9}).

\begin{theorem}\label{thm3.2}
Для нормировочных постоянных $a_n$ и $b_n$
\begin{enumerate}
\item[(a)] имеют место следующие асимптотические соотношения $\left(n \to \infty \right):$
$$%\begin{multline}\label{eq3.4}
a_n \left( q, \alpha, \beta \right) =\dfrac{\pi}{2} \left[ 1+ \dfrac {2 \, s_n \left( q, \alpha, \beta \right)}{\pi \left[n + \delta \left( \alpha, \beta \right)\right]} + r_n \right] \sin^2 \alpha +
$$
$$
+\dfrac{\pi}{2 \left[n + \delta_n(\alpha, \beta)\right]^2} \left[ 1+ \dfrac {2 \, s_n \left( q, \alpha, \beta \right)}{\pi \left[n + \delta \left( \alpha, \beta \right)\right]} + \tilde{r}_n \right]\cos^2 \alpha,
$$%\end{multline}

\begin{multline*}
b_n \left( q, \alpha, \beta \right) =\dfrac{\pi}{2} \left[ 1+ \dfrac {2 \, s_n \left( q, \alpha, \beta \right)}{\pi \left[n + \delta \left( \alpha, \beta \right)\right]} + p_n \right] \sin^2 \beta + \\
+\dfrac{\pi}{2 \left[n + \delta_n(\alpha, \beta)\right]^2} \left[ 1+ \dfrac {2 \, s_n \left( q, \alpha, \beta \right)}{\pi \left[n + \delta \left( \alpha, \beta \right)\right]} + \tilde{p}_n \right]\cos^2 \beta,
\end{multline*}
где
$$%\begin{equation}\label{eq3.5}
s_n = s_n \left( q, \alpha, \beta \right) = -\dfrac {1}{2} \displaystyle \int_{0}^{\pi} \left( \pi - t \right) q\left(t\right) \sin 2 \left[ n + \delta_n\left(\alpha, \beta\right)\right]t dt,
$$%\end{equation}
$r_n = r_n \left( q, \alpha, \beta \right) = O \left( \dfrac {1} {n^2}\right)$ и $\tilde{r}_n=\tilde{r}_n \left( q, \alpha, \beta \right) = O \left( \dfrac {1} {n^2}\right)$ (та же оценка верна для $p_n$ и $\tilde{p}_n$), когда $n \to \infty,$ равномерна по всем $\alpha, \beta \in \left[0, \pi\right]$ и $q \in {BL}^1_\mathbb{R}\left[0, \pi\right].$
\item[(b)] Функция $s,$ определенная формулой
$$%\begin{equation}\label{eq3.6}
s \left(x \right) = \displaystyle \sum_{n=2}^{\infty} \dfrac {s_n}{n + \delta_n \left( \alpha, \beta \right)} \cos \left[n + \delta_n \left( \alpha, \beta \right)\right]x
$$%\end{equation}
абсолютно непрерывна на произвольном отрезке $\left[a,b\right] \subset \left(0, 2 \pi\right),$ т.е. $s \in AC \left(0, 2 \pi\right).$
\end{enumerate}
\end{theorem}

В работе \cite{Harutyunyan_Pahlevanyan:2016} утверждение (b) теоремы \ref{thm3.2} было доказано при условии $\alpha, \beta \in \left(0, \pi \right)$ и в случае $\alpha=\pi,$ $\beta=0.$ Методами, примененными при доказательстве теоремы \ref{thm3.1} можно доказать, что утверждение верно также при $\alpha=\pi,$ $\beta \in \left(0, \pi \right).$

\section{Доказательства теорем \ref{thm2.3} и \ref{thm2.4}}
\label{sec2}

Мы приведем доказательство для Теоремы \ref{thm2.3}. Теорема \ref{thm2.4} может быть доказана аналогично.

\begin{proof}
Для $\left|\lambda\right| \to \infty,$ имеют место следующие асимптотические формулы (\cite{Naimark:1969,Harutyunyan_Hovsepyan:2005,Harutyunyan:2016})
\begin{equation}\label{eq2.3}
\varphi \left(x, \mu, \pi \right) := \varphi_{\pi} \left(x, \mu \right) \equiv \varphi_{\pi} \left(x, \lambda^2 \right)=\dfrac{\sin \lambda x}{\lambda}+O\left(\dfrac{e^{\left|Im \lambda \right|x}}{\left|\lambda \right|^2}\right),
\end{equation}
\begin{equation}\label{eq2.4}
\varphi' \left(x, \mu, \pi \right) := \varphi'_{\pi} \left(x, \mu \right) \equiv \varphi'_{\pi} \left(x, \lambda^2 \right)=\cos \lambda x+O\left(\dfrac{e^{\left|Im \lambda \right|x}}{\left|\lambda \right|}\right),
\end{equation}
$$%
\psi \left(x, \mu, \beta \right) := \psi_{\beta} \left(x, \mu \right) \equiv \psi_{\beta} \left(x, \lambda^2 \right)=\cos \lambda \left( \pi-x \right) \sin \beta + \dfrac{\sin \lambda \left( \pi-x \right)}{\lambda}\cos \beta +
$$
\begin{equation}\label{eq2.5}
+O\left(\dfrac{e^{\left|Im \lambda \right|(\pi-x)}}{\left|\lambda \right|}\right)\sin \beta + O\left(\dfrac{e^{\left|Im \lambda \right|(\pi-x)}}{\left|\lambda \right|^2}\right)\cos \beta,
\end{equation}
$$
\psi' \left(x, \mu, \beta \right) := \psi'_{\beta} \left(x, \mu \right) \equiv \psi'_{\beta} \left(x, \lambda^2 \right)=\left( \lambda \sin \lambda \left( \pi-x \right) + O\left(e^{\left|Im \lambda \right|(\pi-x)} \right) \right)\sin \beta -
$$
\begin{equation}\label{eq2.6}
- \left( \cos \lambda \left( \pi-x \right) + O\left(\dfrac{e^{\left|Im \lambda \right|(\pi-x)}}{\left|\lambda \right|}\right)\right) \cos \beta.
\end{equation}
Из \eqref{eq1.6} и \eqref{eq2.5} следует, что для вронскиана $W_{\pi, \beta}\left(\mu \right)$ мы имеем следующую оценку
\begin{multline}\label{eq2.7}
W_{\pi, \beta}\left(\mu \right) \equiv W_{\pi, \beta}\left(\lambda^2 \right) = -\psi_{\beta} \left(0, \mu \right) = -\cos \lambda \pi \sin \beta - \dfrac{\sin \lambda \pi}{\lambda}\cos \beta + \\
+O\left(\dfrac{e^{\left|Im \lambda \right|\pi}}{\left|\lambda \right|}\right)\sin \beta + O\left(\dfrac{e^{\left|Im \lambda \right|\pi}}{\left|\lambda \right|^2}\right)\cos \beta.
\end{multline}
Обозначим через ${\mathbb{Z}}_{1/6}$ следующую область комплексной плоскости $\mathbb{C}$:
$${\mathbb{Z}}_{1/6}=\left\{\lambda \in \mathbb{C}: \; \left|\lambda-\dfrac{n}{2}\right| \geq \dfrac{1}{6}, \; n \in \mathbb{Z} \right\}.$$
Следующая лемма доказана в \cite{Harutyunyan:2010}, методами которые использовались в \cite{Poschel_Trubowitz:1987}.
\begin{lemma}(\cite{Harutyunyan:2010})\label{lem2.5}
Если $\lambda \in {\mathbb{Z}}_{1/6}$, тогда
\begin{equation}\label{eq2.8}
\left|\sin \pi \lambda \right| \geq \dfrac{1}{7}e^{\left|Im\lambda\right|\pi}, \; \left|\cos \pi \lambda \right| \geq \dfrac{1}{7}e^{\left|Im\lambda\right|\pi}.
\end{equation}
\end{lemma}
Из \eqref{eq2.7} и \eqref{eq2.8} следует что для достаточно большого $\lambda^{*}>0,$  существует константа $C_1>0$ такая что
\begin{equation}\label{eq2.9}
\left|W_{\pi, \beta} \left( \lambda^2 \right)\right| \geq C_1 e^{\left|Im \lambda \right|\pi}\sin \beta, \; \mbox{при} \; \lambda \in {\mathbb{Z}}_{1/6}, \; \left|\lambda \right| > \lambda^{*}.
\end{equation}
Рассмотрим следующую краевую задачу
\begin{equation}\label{eq2.10}
- y'' + q\left( x \right)y = \mu y-f \left( x \right), \;\,\, x \in \left( {0,\pi } \right), \; \mu \in \mathbb{C}, \; f \in L^1 \left[0,\pi\right],
\end{equation}
\begin{equation}\label{eq2.11}
y\left( 0 \right)=0, \; y\left( \pi \right)\cos \beta  + y'\left( \pi \right)\sin \beta  = 0, \; \,\,\beta  \in \left( {0,\pi } \right).
\end{equation}
Хорошо известно, что решение $y\left(x, \mu, f \right)$ краевой задачи \eqref{eq2.10}--\eqref{eq2.11} можно записать в следующей форме (см. например, \cite{Naimark:1969,Levitan_Sargsyan:1988})
\begin{multline}\label{eq2.12}
y \left(x, \mu, f \right) = \dfrac{1}{W_{\pi, \beta}\left( \mu \right)} \psi_{\beta}\left( x, \mu \right) \displaystyle\int_{0}^{x}f\left(t\right)\varphi_{\pi} \left(t, \mu \right)dt+\\
+\dfrac{1}{W_{\pi, \beta}\left( \mu \right)}\varphi_{\pi} \left(x, \mu \right) \displaystyle\int_{x}^{\pi}f\left(t\right)\psi_{\beta} \left(t, \mu \right)dt.
\end{multline}
Так как $\varphi,$ $\psi$ и $W_{\pi, \beta}$ являются целыми функциями от $\mu$, то $y \left(x, \mu, f \right)$ является мероморфной функцией от $\mu,$ с полюсами в нулях функции $W_{\pi, \beta}$ или, что то же самое, в собственных значениях $\mu_{n}, n=0,1,2,\dots.$ Поскольку $\dot{W}_{\pi, \beta}\left(\mu_n\right) \equiv \dfrac{d}{d\mu}W_{\pi, \beta} \left(\mu_n\right)=\beta_n a_n$ (см. \cite[Лемма 1.1.1]{Yurko:2007}), то используя \eqref{eq1.8}, мы получаем вычет
\begin{equation}\label{eq2.13}
\mathop {\operatorname{Res}}\limits_{\mu  = {\mu _n}} {y \left(x, \mu, f \right)}
=\dfrac{1}{a_n}\varphi_\pi \left( x, \mu_n \right)\displaystyle\int_{0}^{\pi}f\left( t \right) \varphi_\pi \left(t, \mu_n \right)dt.
\end{equation}
Из \eqref{eq2.3}, \eqref{eq2.5}, \eqref{eq2.9} и \eqref{eq2.12} следует,
что существуют положительные числа $C,$ $C_2,$ $C_3,$ $C_4$ такие что при $\lambda \in {\mathbb{Z}}_{1/6},$ $\left|\lambda \right| > \lambda^{*}$ имеет место следующая оценка
\begin{multline}\label{eq2.14}
\left| y \left(x, \lambda^2, f \right) \right| \leq \dfrac{\left| \psi_{\beta} \left(x, \lambda^2 \right) \right|\mathop {\operatorname{max}}\limits_{t \in \left[0, x \right]}\left|\varphi_{\pi} \left(t, \lambda^2 \right) \right|\displaystyle\int_{0}^{x}\left|f\left( t \right)\right|dt}{C_1 e^{\left| Im \lambda \right|\pi}\sin \beta}+\\
+\dfrac{\left|\varphi_{\pi} \left(x, \lambda^2 \right) \right|\mathop {\operatorname{max}}\limits_{t \in \left[x, \pi \right]}\left|\psi_{\beta} \left(t, \lambda^2 \right) \right| \displaystyle\int_{x}^{\pi}\left|f\left( t \right)\right|dt}{C_1 e^{\left| Im \lambda \right|\pi} \sin \beta} \leq \\
\leq \dfrac{e^{\left| Im \lambda \right|(\pi-x)} \left(\sin \beta+\dfrac{ \left| \cos \beta \right| }{\left| \lambda \right|}+C_3\dfrac{ \sin \beta }{\left| \lambda \right|}+C_4\dfrac{ \left| \cos \beta \right| }{\left| \lambda \right|^2} \right)}{C_1 e^{\left| Im \lambda \right|\pi} \sin \beta} \times \\
\times e^{\left| Im \lambda \right|x} \left(\dfrac{1}{\left| \lambda \right|}+C_2\dfrac{1}{\left| \lambda \right|^2} \right)\displaystyle\int_{0}^{\pi}\left|f\left( t \right)\right|dt \leq \\
\leq \dfrac{1}{C_1} \displaystyle\int_{0}^{\pi}\left|f\left( t \right)\right|dt \left( \dfrac{1}{\left| \lambda \right|}+O\left( \dfrac{1}{\left| \lambda \right|^2} \right) \right) \leq \dfrac{C}{\left| \lambda \right|}.
\end{multline}
Рассмотрим теперь функцию $f \in AC\left[0, \pi \right].$ Используя тот факт, что $\varphi_{\pi} \left(x, \mu \right)$ и $\psi_{\beta} \left(x, \mu \right)$ являются решениями \eqref{eq1.1}, мы можем переписать представление \eqref{eq2.12} для $y \left(x, \mu, f \right)$ в следующем виде (сравните с \cite{Yurko:2007}):
\begin{equation}\label{eq2.15}
y \left(x, \mu, f \right) = \dfrac{f\left(x\right)}{\mu}+f\left(0\right) \dfrac{\psi_{\beta}\left(x, \mu \right)}{\mu W_{\pi, \beta}\left(\mu\right)}+\dfrac{Z_1\left(x, \mu, \pi, \beta, f' \right)} {\mu}+\dfrac{Z_2\left(x, \mu, \pi, \beta \right)}{\mu},
\end{equation}
где
\begin{equation}\label{eq2.16}
Z_1\left(x, \mu, \pi, \beta, f' \right)=\dfrac{\psi_{\beta}\left( x, \mu \right) \displaystyle\int_{0}^{x} f'\left(t\right) \varphi'_{\pi} \left(t, \mu \right)dt +\varphi_{\pi}\left( x, \mu \right) \displaystyle\int_{x}^{\pi}f'\left(t\right)\psi'_{\beta} \left(t, \mu \right)dt}{W_{\pi, \beta}\left( \mu \right)},
\end{equation}
\begin{multline}\label{eq2.17}
Z_2\left(x, \mu, \pi, \beta \right)=-f\left(\pi\right)\psi'_{\beta} \left(\pi, \mu \right) \dfrac{\varphi_{\pi}\left( x, \mu \right)}{W_{\pi, \beta}\left( \mu \right)}+y \left(x, \mu, qf \right)= \\
=f\left(\pi\right)\cos \beta \dfrac{\varphi_{\pi}\left( x, \mu \right)}{W_{\pi, \beta}\left( \mu \right)}+y \left(x, \mu, qf \right).
\end{multline}
Покажем, что
\begin{equation}\label{eq2.18}
\mathop {\lim }\limits_{\left| \lambda  \right| \to \infty \hfill \atop \lambda  \in {\mathbb{Z}}_{1/6}} \hfill \mathop {\operatorname{max}}\limits_{x \in \left[0, \pi\right]}\left|Z_1\left(x, \mu, \pi, \beta, f' \right)\right|=0.
\end{equation}
Сначала предположим, что $f'$ абсолютно непрерывная функция на $\left[0, \pi\right].$ Тогда существует $f'' \in L^1\left[ {0,\pi } \right]$ и \eqref{eq2.16} можно записать в следующем виде
\begin{multline*}
Z_1\left(x, \mu, \pi, \beta, f' \right)=\dfrac{\varphi_{\pi}\left(x, \mu \right)}{W_{\pi, \beta}\left( \mu \right)} f'\left(\pi \right) \sin \beta-\\
-\dfrac{\psi_{\beta}\left( x, \mu \right) \displaystyle\int_{0}^{x} f''\left(t\right) \varphi_{\pi} \left(t, \mu \right)dt +\varphi_{\pi}\left( x, \mu \right) \displaystyle\int_{x}^{\pi}f''\left(t\right)\psi_{\beta} \left(t, \mu \right)dt}{W_{\pi, \beta}\left( \mu \right)}.
\end{multline*}
В силу \eqref{eq2.3}--\eqref{eq2.6} и \eqref{eq2.9} мы получаем, что существует число $C>0,$ такое что
\begin{equation*}
\mathop {\operatorname{max}}\limits_{x \in \left[0, \pi\right]}\left|Z_1\left(x, \mu, \pi, \beta, f' \right)\right| \leq \dfrac{C}{\left| \lambda \right|}, \; \mbox{при} \; \lambda \in {\mathbb{Z}}_{1/6}, \; \left|\lambda \right| > \lambda^{*}.
\end{equation*}
Отсюда следует \eqref{eq2.18} в случае $f' \in AC\left[0,\pi\right].$

Теперь обратимся к общему случаю $g := f' \in L^1 \left[0, \pi \right].$ Зафиксируем $\epsilon>0$ и выберем абсолютно непрерывную функцию $g_\epsilon,$ так что
\begin{equation*}
\displaystyle\int_{0}^{\pi} \left|g\left(t\right)-{g_\epsilon}\left(t\right)\right|dt<\dfrac{C_1 \sin \beta }{16} \, \epsilon.
\end{equation*}
Тогда, согласно \eqref{eq2.3}--\eqref{eq2.6}, \eqref{eq2.9} и \eqref{eq2.16} для $\lambda \in {\mathbb{Z}}_{1/6}, \; \left|\lambda \right| > \lambda^{*},$ мы имеем
\begin{multline*}
\mathop {\operatorname{max}}\limits_{x \in \left[0, \pi\right]}\left|Z_1\left(x, \mu, \pi, \beta, g \right)\right| \leq \mathop {\operatorname{max}}\limits_{x \in \left[0, \pi\right]}\left|Z_1\left(x, \mu, \pi, \beta, {g_\epsilon}\right)\right|+\mathop {\operatorname{max}}\limits_{x \in \left[0, \pi\right]}\left|Z_1\left(x, \mu, \pi, \beta, g-{g_\epsilon}\right)\right| \leq \\
\leq \dfrac{C(\epsilon)}{\left| \lambda \right|}+\dfrac{C_1 \sin \beta }{16} \, \epsilon \mathop {\operatorname{max}}\limits_{x \in \left[0, \pi\right]}\left(\dfrac{\left|\psi_{\beta}\left(x, \mu \right)\right|\mathop {\operatorname{max}}\limits_{t \in \left[0, x \right]}\left|\varphi'_{\pi}\left(t, \mu\right)\right|+\left|\varphi_{\pi}\left(x, \mu \right)\right|\mathop {\operatorname{max}}\limits_{t \in \left[0, x \right]}\left|\psi'_{\beta}\left(t, \mu\right)\right|}{C_1 e^{\left|Im \lambda \right|\pi}\sin \beta}\right)\leq \\
\leq \dfrac{C(\epsilon)}{\left| \lambda \right|}+\dfrac{C_1 \sin \beta }{16} \, \epsilon \mathop {\operatorname{max}}\limits_{x \in \left[0, \pi\right]}\left(\dfrac{8 e^{\left|Im \lambda \right|\pi}}{C_1 e^{\left|Im \lambda \right|\pi}\sin \beta}\right) \leq \dfrac{C(\epsilon)}{\left| \lambda \right|}+\dfrac{\epsilon}{2}.
\end{multline*}
Легко видеть, что если мы выберем $\lambda_\epsilon^{*}=\dfrac{2C\left(\epsilon\right)}{\epsilon},$ тогда для $\lambda \in {\mathbb{Z}}_{1/6}$ и $\left|\lambda \right| > \lambda_\epsilon^{*}$ мы имеем $\mathop {\operatorname{max}}\limits_{x \in \left[0, \pi\right]}\left|Z_1\left(x, \mu, \pi, \beta\right)\right| \leq \epsilon.$  В силу произвольности $\epsilon>0,$ мы приходим к \eqref{eq2.18}.
Теперь оценим $Z_2\left(x, \mu, \pi, \beta\right)$ (см. \eqref{eq2.17}). Поскольку $qf \in L^1\left[ {0,\pi } \right],$ то оценки в \eqref{eq2.14} верны также для $y\left( x, \mu, qf \right).$ Используя \eqref{eq2.3}, \eqref{eq2.9}, \eqref{eq2.14} и тот факт, что $\sin \beta \neq 0$ мы получаем следующие оценки (при $\lambda \in {\mathbb{Z}}_{1/6}, \; \left|\lambda \right| > \lambda^{*}$):
\begin{multline}\label{eq2.19}
\mathop {\operatorname{max}}\limits_{x \in \left[0, \pi\right]}\left|Z_2\left(x, \mu, \pi, \beta\right)\right| \leq \mathop {\operatorname{max}}\limits_{x \in \left[0, \pi\right]} \left|f\left(\pi\right)\cos \beta \dfrac{\varphi_{\pi}\left( x, \mu \right)}{W_{\pi, \beta}\left( \mu \right)}\right| + \mathop {\operatorname{max}}\limits_{x \in \left[0, \pi\right]}\left|y \left(x, \mu, qf \right)\right| \leq \\
\leq \left|f\left(\pi\right)\cos \beta \dfrac{C_5 e^{\left|Im \lambda\right|\pi}}{\left|\lambda\right|C_1 e^{\left|Im \lambda\right|\pi}\sin \beta}\right| +\dfrac{C_6}{\left|\lambda\right|} \leq \dfrac{C_5}{C_1}\dfrac{\left|f\left(\pi\right)\cot \beta\right|} {\left|\lambda\right|} +\dfrac{C_6}{\left|\lambda\right|} \leq \dfrac{C_7}{\left|\lambda\right|},
\end{multline}
где $C_5-C_7$ положительные числа. Рассмотрим следующий контурный интеграл
$$%\begin{equation}\label{eq2.20}
I_{N}\left(x\right)=\dfrac{1}{2\pi i}\oint_{\Gamma_{N}}y \left(x, \mu, f \right)d \mu,
$$%\end{equation}
где $\Gamma_{N}=\left\{\mu : \left|\mu\right|=\left(N+\dfrac{3}{4}\right)^2\right\}$ (с обходом против часовой стрелки).
С одной стороны, используя теорему Коши о вычетах (см. \cite{Shabat:1985}), из \eqref{eq2.13} мы получаем
\begin{equation}\label{eq2.21}
I_{N}\left(x\right)= \sum\limits_{n = 0}^{N} \dfrac{1}{a_n}\displaystyle\int_{0}^{\pi}f\left( t \right){\varphi_{\pi} \left(t, \mu_n \right)}dt \varphi_{\pi} \left(x, \mu_n \right).
\end{equation}
С другой стороны, из \eqref{eq2.15}, \eqref{eq2.18} и \eqref{eq2.19} имеем, что
\begin{equation}\label{eq2.22}
I_{N}\left(x\right)= f\left(x\right)+f(0)\dfrac{1}{2\pi i}\oint_{\Gamma_{N}} \dfrac{\psi_{\beta}\left(x, \mu \right)}{\mu W_{\pi, \beta}\left(\mu\right)}d \mu+{\epsilon_{N}}\left(x\right),
\end{equation}
где ${\epsilon_{N}}\left(x\right),$ согласно \eqref{eq2.18} и \eqref{eq2.19}, равномерно сходится к $0:$
$$ \mathop {\lim}\limits_{N \to \infty}\mathop {\operatorname{max}}\limits_{x \in \left[0, \pi\right]}\left|{\epsilon_{N}}\left(x\right)\right|=0.$$
Без потери общности, будем считать что $\mu=0$ не является собственным значением задачи $L\left(q,\pi,\beta\right).$ В самом деле, из чистой дискретности спектра следует, что существует число $c$ такое, что числа $\mu_n +c \neq 0,$ $n=0,1,2,\dots,$ являются собственными значениями задачи $L\left(q+c, \pi, \beta \right)$ с теми же собственными функциями $\varphi_n$ и нормировочными постоянными $a_n,$ что и у задачи $L\left(q, \pi, \beta \right).$ Тогда функция $\dfrac{\psi_{\beta}\left(x, \mu \right)}{\mu W_{\pi, \beta}\left(\mu\right)}$ имеет полюсы только первого порядка и используя теорему Коши о вычетах мы можем легко вычислить, что
$$
\phi_N \left(x \right) :=\dfrac{1}{2\pi i}\oint_{\Gamma_{N}} \dfrac{\psi_{\beta}\left(x, \mu \right)}{\mu W_{\pi, \beta}\left(\mu\right)}d \mu=\mathop {\operatorname{Res}}\limits_{\mu=0}\dfrac{\psi_{\beta}\left(x, \mu \right)}{\mu W_{\pi, \beta}\left(\mu\right)} +\sum\limits_{n = 0}^{N}\mathop {\operatorname{Res}}\limits_{\mu=\mu_n}\dfrac{\psi_{\beta}\left(x, \mu \right)}{\mu W_{\pi, \beta}\left(\mu\right)}=
$$
$$
=\dfrac{\psi_{\beta}\left(x, 0 \right)}{W_{\pi, \beta}\left(0\right)}+\sum\limits_{n = 0}^{N} \dfrac{\psi_{\beta}\left(x, \mu_n \right)}{\mu_n \dot{W}_{\pi, \beta}\left(\mu_n\right)}= \dfrac{\psi_{\beta}\left(x, 0 \right)}{W_{\pi, \beta}\left(0\right)}+\sum\limits_{n = 0}^{N} \dfrac{\beta_n \varphi_{\pi}\left(x, \mu_n \right)}{\mu_n \beta_n a_n}=
$$
\begin{equation}\label{eq2.23}=\dfrac{\psi_{\beta}\left(x, 0 \right)}{W_{\pi, \beta}\left(0\right)}+\sum\limits_{n = 0}^{N} \dfrac{1}{\mu_n a_n} {\varphi_n \left(x \right)}.
\end{equation}
Теперь покажем, что последовательность $\phi_N \left(x \right)$ сходится к $0$ (при $N \to \infty$) равномерно на сегменте $\left[a, \pi \right],$ для произвольного $a \in \left(0, \pi \right).$

Так как $\varphi_n\left(x\right)=\dfrac{\sin \left(n+\frac{1}{2}\right)x}{n+\frac{1}{2}}+O\left(\dfrac{1}{n^2}\right)$ равномерно на $\left[0, \pi \right]$ (см. \eqref{eq2.3}), $\mu_n=\mu_n \left(q,\pi,\beta\right)=\left(n+\dfrac{1}{2}\right)^2+O\left(1\right)$ (см. \cite[Теорема 1 на стр. 286 и асимптотические оценки для $\delta_n \left(\pi, \beta \right)$ на стр. 292]{Harutyunyan:2008}) и $a_n=a_n\left(q,\pi,\beta\right)=\dfrac{\pi}{2\left(n+\frac{1}{2}\right)^2} \left(1+o\left(\dfrac{1}{n}\right)\right)$ (см. \cite[Теорема 1.1 ]{Harutyunyan_Pahlevanyan:2016}, стр. 9-10), тогда $\phi_N \left(x \right)$ (см. \eqref{eq2.23}) можно записать в следующем виде:
\begin{equation*}
\phi_N \left(x \right)=\dfrac{\psi_{\beta}\left(x, 0 \right)}{W_{\pi, \beta}\left(0\right)}+\dfrac{2}{\pi}\sum\limits_{n = 0}^N \dfrac{\sin \left(n+\frac{1}{2}\right)x}{n+\frac{1}{2}}+\sum\limits_{n = 0}^N q_n \left(x\right),
\end{equation*}
где $q_n \left(x \right)=O\left(\dfrac{1}{n^2}\right)$ равномерно на $\left[0, \pi \right].$

Поскольку $\displaystyle\sum\limits_{n = 0}^\infty \dfrac{\sin \left(n+\frac{1}{2}\right)x}{n+\frac{1}{2}}=\dfrac{\pi}{2},$ $0 < x < 2\pi$ (см. например, \cite[формула (37) на стр. 578]{Bronshtein_Semendyayev:1998}), тогда последовательность $\phi_N \left(x \right)$ сходится к непрерывной функции $\phi \left(x \right)$ (при $N \to \infty$) равномерно на сегменте $\left[a, \pi \right],$ для произвольного $a \in \left(0, \pi \right).$

Теперь, чтобы доказать, что $\phi \left(x \right) \equiv 0, \, x \in \left(0, \pi \right],$ достаточно показать, что $\phi=0$ п.в.

Сделав некоторые вычисления, получим:
\begin{equation}\label{eq2.24}
\int\limits_0^\pi \phi \left(x \right) \varphi_m \left(x \right)dx = \dfrac{1}{W_{\pi, \beta}\left(0\right)} \int\limits_0^\pi \psi_{\beta}\left(x, 0 \right) \varphi_m \left(x\right)dx+\dfrac{1}{\mu_m}, \, m=0,1,2,\dots
\end{equation}
\begin{multline}\label{eq2.25}
\mu_m \int\limits_0^\pi \psi_{\beta}\left(x, 0 \right) \varphi_m \left(x\right)dx=\int\limits_0^\pi \left(\varphi_m \left(x\right)\psi''_{\beta}\left(x, 0 \right)-\varphi''_m \left(x\right) \psi_{\beta}\left(x, 0 \right)\right)dx=\\
=\left.\left(\varphi_m \left(x\right)\psi'_{\beta}\left(x, 0 \right)-\varphi'_m \left(x\right) \psi_{\beta}\left(x, 0 \right)\right)\right|_{0}^{\pi}=\psi_{\beta}\left(0, 0 \right)=-W_{\pi, \beta}\left(0\right).
\end{multline}

Из \eqref{eq2.24} и \eqref{eq2.25} следует
\begin{equation*}
\int\limits_0^\pi \phi \left(x \right) \varphi_m \left(x \right)dx =0, \, m=0,1,2,\dots.
\end{equation*}

Поскольку система собственных функций $\left\{\varphi_{m} \left(x \right)\right\}_{m=0}^\infty$ 
краевой задачи $L(q, \pi, \beta)$ является полной и ортогональной в $L^2 \left(0, \pi \right)$, то $\phi=0$ п.в.
Сравнивая этот результат с \eqref{eq2.21}, \eqref{eq2.22} и переходя 
к пределу при $N \to \infty $ в \eqref{eq2.22}, мы приходим к \eqref{eq2.1}. Теорема \ref{thm2.3} доказана.
\end{proof}

\begin{remark}\label{rem2.6}
Хорошо известно, что одно из доказательств теоремы \ref{thm2.2} основывается на так называемой теореме о равномерной 
равносходимости, которая утверждает, что разложение по собственным функциям задачи $L \left(q, \alpha, \beta \right),$ $\alpha, \beta \in \left(0,\pi \right)$ эквивалентно разложению по собственным функциям задачи

\noindent$L \left(0, \dfrac{\pi}{2}, \dfrac{\pi}{2} \right),$ т.е., $\left\{\cos nx\right\}_{n \geq 0}$ (см. \cite{Coddington_Levinson:1955,Levitan_Sargsyan:1970,Levitan_Sargsyan:1988}). Далее можно применить теорему Дирихле-Жордана (см. \cite[стр. 121--122]{Bari:1961}) и Теорема \ref{thm2.2} будет доказана. Тот же подход не может быть применен в нашем случае, а именно: нетрудно установить что разложение по собственным функциям задачи $L \left(q, \pi, \beta \right),$ $\beta \in \left(0, \pi \right)$ будет эквивалентно разложению по $\sin \left(n+\dfrac{1}{2}\right)x,$ $n=0,1,2,\dots,$ т.е., по собственным функциям задачи $L \left(0, \pi, \dfrac{\pi}{2} \right)$ (см. \cite[замечание на стр. 304]{Coddington_Levinson:1955} и \cite[замечание на стр. 71]{Levitan_Sargsyan:1970}). С другой стороны, насколько нам известно, нет аналога теоремы Дирихле-Жордана для разложения по системе функций $\left\{\sin \left(n+\dfrac{1}{2}\right)x\right\}_{n \geq 0}$ (по этому поводу см. \cite[Теорема 2.6]{Iserles_Norsett:2008}).
\end{remark}

\section{Асимптотика собственных значений}
\label{sec3}

Стоит заметить что приведенное в работе \cite{Harutyunyan:2016} доказательство утверждения (b) теоремы \ref{thm3.1} для случая $\alpha, \beta \in \left(0,\pi \right)$ не проходит для случая $\alpha=\pi,$ $\beta \in \left(0,\pi \right).$ Ниже, используя теоремы \ref{thm2.3} и \ref{thm2.4}, мы разберем этот случай.
Обозначим $\sigma\left(x\right)=\displaystyle\int_{0}^{x} q\left( t \right)dt$ и запишем $l_n\left(q, \pi, \beta \right)$ (см. \eqref{eq3.2}) в следующей форме:
$$%\begin{equation}\label{eq3.7}
l_n\left(q, \pi, \beta\right)=-\dfrac{\sigma \left(\pi\right)\cos 2\pi \delta_n\left(\pi, \beta\right)}{2 \pi \left(n+\delta_n\left(\pi, \beta\right)\right)}-\dfrac{1}{2\pi}\int_{0}^{2\pi}\sigma_1 \left(x\right)\sin\left(n+\delta_n\left(\pi, \beta\right)\right)xdx,
$$%\end{equation}
где $\sigma_1 \left( x \right)\equiv \sigma \left( \dfrac{x}{2} \right)$ абсолютно непрерывная функция на $\left[0, 2\pi\right].$

Заметим что единственность решения $\delta_n \left(\alpha, \beta \right)$ уравнения \eqref{eq1.10} при $\alpha=\pi$ и $\beta \in \left[0, \pi \right)$ можно доказать исходя из того, что $\arccos$ является убывающей функцией.

Из \eqref{eq1.10} легко видеть (подробности см. \cite{Harutyunyan:2008}), что для $\beta \in \left(0,\pi\right)$ мы имеем \begin{equation}\label{eq3.8}
\delta_n \left(\pi, \beta \right)=\dfrac{1}{2}+\dfrac{\cot\beta}{\pi\left(n+\frac{1}{2}\right)}+ O\left(\dfrac{1}{n^2}\right)\,{\cot\beta} =\dfrac{1}{2}+O\left(\dfrac{1}{n}\right),
\end{equation}
и следовательно,
\begin{equation}\label{eq3.9}
\cos 2 \pi \delta_n\left(\pi, \beta\right)=-1+d_n, \; \sin 2 \pi \delta_n\left(\pi, \beta\right)=e_n,
\end{equation}
где $d_n=O\left(\dfrac{1}{n^2}\right), \; e_n=O\left(\dfrac{1}{n}\right).$

Поэтому, $l\left(x, \beta \right)$ (см. \eqref{eq3.3}) можно представить в виде суммы трех функций
$$%\begin{equation}\label{eq3.10}
l\left(x, \beta\right)=l_1\left(x, \beta\right)+l_2\left(x, \beta\right)+l_3\left(x, \beta\right),
$$%\end{equation}
где
$$%\begin{equation}\label{eq3.11}
l_1\left(x, \beta\right)=\dfrac{\sigma\left(\pi\right)}{2 \pi}\displaystyle\sum_{n=2}^{\infty} \dfrac{\sin \left(n+\delta_n\left(\pi, \beta\right)\right)x}{\left(n+\delta_n\left(\pi, \beta\right)\right)},
$$%\end{equation}
\begin{equation}\label{eq3.12}
l_2\left(x, \beta\right)=-\dfrac{\sigma\left(\pi\right)}{2 \pi}\displaystyle\sum_{n=2}^{\infty} {d_n}\dfrac{\sin \left(n+\delta_n\left(\pi, \beta\right)\right)x}{\left(n+\delta_n\left(\pi, \beta\right)\right)},
\end{equation}
\begin{equation}\label{eq3.13}
l_3\left(x, \beta\right)=-\dfrac{1}{2 \pi} \displaystyle\sum_{n=2}^{\infty} f_n \sin \left(n+\delta_n\left(\pi, \beta\right)\right)x,
\end{equation}
и $f_n=\displaystyle\int_{0}^{2\pi} \sigma_1 \left(t\right)\sin\left(n+\delta_n\left(\pi, \beta\right)\right)tdt.$ \\
Поскольку $f_n=\displaystyle\int_{0}^{\pi} \sigma_1 \left(t\right)\sin\left(n+\delta_n\left(\pi, \beta\right)\right)tdt+\displaystyle\int_{\pi}^{2\pi} \sigma_1 \left(t\right)\sin\left(n+\delta_n\left(\pi, \beta\right)\right)tdt$ и \\
\begin{multline*}
\int_{\pi}^{2\pi} \sigma_1 \left(t\right)\sin\left(n+\delta_n\left(\pi, \beta\right)\right)tdt=
\int_{-2\pi}^{-\pi} -\sigma_1 \left(-t\right)\sin\left(n+\delta_n\left(\pi, \beta\right)\right)tdt=\\
=\int_{0}^{\pi} -\sigma_1 \left(2\pi-t\right)\sin\left(n+\delta_n\left(\pi, \beta\right)\right)\left(t-2\pi\right)dt=\\
=\int_{0}^{\pi} \sigma_1 \left(2\pi-t\right)\left(\left(1-d_n\right)\sin\left(n+\delta_n\left(\pi, \beta\right)\right)t+e_n \cos\left(n+\delta_n\left(\pi, \beta\right)\right)t \right)dt=\\
=\int_{0}^{\pi} \sigma \left(\pi-\dfrac{t}{2}\right)\left(\left(1-d_n\right)\sin\left(n+\delta_n\left(\pi, \beta\right)\right)t+e_n \cos\left(n+\delta_n\left(\pi, \beta\right)\right)t \right)dt,
\end{multline*}
то
\begin{multline}\label{eq3.14}
f_n=\int_{0}^{\pi} \left(\sigma \left(\dfrac{t}{2}\right)+\sigma \left(\pi-\dfrac{t}{2}\right)\right)\sin\left(n+\delta_n\left(\pi, \beta\right)\right)tdt-\\
-d_n \int_{0}^{\pi} \sigma \left(\pi-\dfrac{t}{2}\right)\sin\left(n+\delta_n\left(\pi, \beta\right)\right)tdt+
e_n \int_{0}^{\pi} \sigma \left(\pi-\dfrac{t}{2}\right)\cos\left(n+\delta_n\left(\pi, \beta\right)\right)t.
\end{multline}
Следует отметить, что $\delta_n\left(\alpha, \beta\right)$ определена только для $n \geq 2,$ поэтому мы запишем $\lambda_0 \left(0,\pi,\beta\right),$ $\lambda_1 \left(0,\pi,\beta\right)$ и $\lambda_n \left(0,\pi,\beta\right)=n+\delta_n\left(\pi, \beta\right)$ для всех $n \geq 2.$
Учитывая, что система функций
$$\left\{\varphi_{n} \left(x \right)\right\}_{n=0}^\infty=
%\left\{\dfrac{\sin \lambda_n \left(0,\pi,\beta\right)x}{\lambda_n \left(0,\pi,\beta\right)}\right\}_{n=0}^\infty
 \left\{\dfrac{\sin \lambda_n \left(0,\pi,\beta\right)x}{\lambda_n \left(0,\pi,\beta\right)}\right\}_{n=0}^1 \cup \left\{\dfrac{\sin \left(n+\delta_n\left(\pi, \beta \right) \right)x}{n+\delta_n\left(\pi, \beta \right)}\right\}_{n=2}^\infty$$
является системой собственных функций задачи $L\left(0, \pi, \beta\right)$ и применяя Теорему \ref{thm2.3}, получаем
\begin{multline}\label{eq3.15}
\sigma \left(\dfrac{x}{2}\right)+\sigma \left(\pi-\dfrac{x}{2}\right) =\sigma_{2}\left(x\right)+\\
+\sum\limits_{n = 2}^\infty \dfrac{\displaystyle\int_{0}^{\pi}\left(\sigma \left(\dfrac{t}{2}\right)+\sigma \left(\pi-\dfrac{t}{2}\right)\right)\sin \left(n+\delta_n\left(\pi, \beta \right) \right)tdt}{\displaystyle\int_{0}^{\pi}\sin^{2}\left(n+\delta_n\left(\pi, \beta \right) \right) tdt} \sin \left(n+\delta_n\left(\pi, \beta \right) \right)x
\end{multline}
где ряд сходится равномерно на произвольном отрезке $\left[a, \pi \right] \subset \left(0, \pi \right]$ и
\begin{equation*}
\sigma_{2}\left(x\right):= \sum\limits_{n = 0}^1 \dfrac{\displaystyle\int_{0}^{\pi}\left(\sigma \left(\dfrac{t}{2}\right)+\sigma \left(\pi-\dfrac{t}{2}\right)\right)\sin \lambda_n \left(0,\pi,\beta\right)tdt}{\displaystyle\int_{0}^{\pi}\sin^{2}\lambda_n \left(0,\pi,\beta\right)tdt}\sin \lambda_n \left(0,\pi,\beta\right)x
\end{equation*}
Используя \eqref{eq3.8} и \eqref{eq3.9}, мы вычисляем
\begin{multline}\label{eq3.16}
\int\limits_0^\pi {\sin^2 \left(n+{\delta_n}\left( \pi, \beta \right)\right)t}dt=\\
=\dfrac{\pi}{2}-\dfrac{\sin 2 \pi \left(n+\delta_n (\pi, \beta)\right)}{4\left(n+\delta_n (\pi, \beta)\right)}=\dfrac{\pi}{2} - \dfrac{e_n}{4\left(n+\delta_n (\pi, \beta)\right)}.
\end{multline}
Из \eqref{eq3.16}, легко видеть, что
$$
\dfrac{1}{\int\limits_0^\pi {\sin^2 \left(n+{\delta_n}\left( \pi, \beta \right) \right)t}dt} = \dfrac{2}{\pi}+{g_n},
\quad\text{где}\quad g_n=\dfrac{2 e_n}{\pi\left(2\pi\left(n+\delta_n\left(\pi, \beta \right) \right)-e_n \right)}=O\left(\dfrac{1}{n^2}\right).
$$
Теперь мы можем записать \eqref{eq3.15} в форме
\begin{multline}\label{eq3.17}
\sum\limits_{n = 2}^\infty \dfrac{2}{\pi}\displaystyle\int_{0}^{\pi}\left(\sigma \left(\dfrac{t}{2}\right)+\sigma \left(\pi-\dfrac{t}{2}\right)\right)\sin \left(n+\delta_n\left(\pi, \beta \right) \right)tdt \sin \left(n+\delta_n\left(\pi, \beta \right) \right)x=\\
=-\sum\limits_{n = 2}^\infty {g_n}\displaystyle\int_{0}^{\pi}\left(\sigma \left(\dfrac{t}{2}\right)+\sigma \left(\pi-\dfrac{t}{2}\right)\right)\sin \left(n+\delta_n\left(\pi, \beta \right) \right)tdt \sin \left(n+\delta_n\left(\pi, \beta \right) \right)x+\\
+\sigma \left(\dfrac{x}{2}\right)+\sigma \left(\pi-\dfrac{x}{2}\right) -\sigma_2\left(x\right),
\end{multline}
где ряды сходятся равномерно на произвольном отрезке $\left[a, \pi \right] \subset \left(0, \pi \right].$

Из \eqref{eq3.13}, \eqref{eq3.14}, \eqref{eq3.17} следует что для произвольного $x \in \left(0, \pi \right]$
\begin{multline}\label{eq3.18}
l_3\left(x, \beta\right) = \dfrac{1}{4}\left(-\sigma \left(\dfrac{x}{2}\right)-\sigma \left(\pi-\dfrac{x}{2}\right)+\sigma_2\left(x\right)\right)+\\
+\displaystyle\sum_{n=2}^{\infty} \dfrac{g_n}{4} \displaystyle\int_{0}^{\pi}\left(\sigma \left(\dfrac{t}{2}\right)+\sigma \left(\pi-\dfrac{t}{2}\right)\right)\sin \left(n+\delta_n\left(\pi, \beta \right) \right)tdt \sin \left(n+\delta_n\left(\pi, \beta \right) \right)x+\\
+\dfrac{1}{2\pi}\displaystyle\sum_{n=2}^{\infty} {d_n} \int_{0}^{\pi} \sigma \left(\pi-\dfrac{t}{2}\right)\sin\left(n+\delta_n\left(\pi, \beta\right)\right)tdt \sin \left(n+\delta_n\left(\pi, \beta\right)\right)x - \\
-\dfrac{1}{2\pi}\displaystyle\sum_{n=2}^{\infty} {e_n} \int_{0}^{\pi} \sigma \left(\pi-\dfrac{t}{2}\right)\cos\left(n+\delta_n\left(\pi, \beta\right)\right)tdt \sin \left(n+\delta_n\left(\pi, \beta\right)\right)x.
\end{multline}
Поскольку $d_n=O\left(\dfrac{1}{n^2}\right),$ $e_n \displaystyle\int_{0}^{\pi} \sigma \left(\pi-\dfrac{t}{2}\right)\cos\left(n+\delta_n\left(\pi, \beta\right)\right)tdt=O\left(\dfrac{1}{n^2}\right),$ $g_n=O\left(\dfrac{1}{n^2}\right),$ тогда $l_3 \in AC\left(0, \pi\right].$ С другой стороны, так как (см. \eqref{eq3.13} и \eqref{eq3.9})
\begin{multline*}
l_3\left(2\pi-x, \beta\right)=l_3\left(x, \beta\right)+\dfrac{1}{2 \pi} \displaystyle\sum_{n=2}^{\infty} d_n f_n \sin \left(n+\delta_n\left(\pi, \beta\right)\right)x-\\
-\dfrac{1}{2 \pi} \displaystyle\sum_{n=2}^{\infty} e_n f_n \cos \left(n+\delta_n\left(\pi, \beta\right)\right)x,
\end{multline*}
то $l_3 \in AC\left[\pi, 2\pi\right)$ и следовательно $l_3 \in AC\left(0, 2\pi\right)$. Поскольку $\displaystyle\sum_{n=2}^{\infty} \dfrac{\sin \left(n+\delta_n\left(\pi, \beta\right)\right)x}{\left(n+\delta_n\left(\pi, \beta\right)\right)}$ абсолютно непрерывная функция на $\left(0, 2 \pi\right)$ (см. \cite{Harutyunyan:2010,Harutyunyan:2016}), тогда $l_1 \in AC\left(0, 2\pi\right)$.
Поскольку $d_n=O\left(\dfrac{1}{n^2}\right),$ то ряд в \eqref{eq3.12} и его первая производная сходятся абсолютно и равномерно на $\left[0, 2\pi\right]$ и, следовательно $l_2 \in AC\left[0, 2\pi\right].$ Утверждение (b) теоремы \ref{thm3.1} при $\alpha=\pi,$ $\beta \in \left(0,\pi \right)$ доказано.
\hfill$\Box$

\subsection*{Благодарность} Автор выражает благодарность профессору Т. Н. Арутюняну за постановку задачи и внимание к работе.

\vskip4mm

\noindent{\bf Abstract.}\, Uniform convergence of the expansion of an absolutely continuous function
for eigenfunctions of the Sturm-Liouville problem $-y'' + q \left( x \right) y = \mu y,$ $y \left(0\right)=0,$
$y\left( \pi \right)\cos \beta + y'\left( \pi \right)\sin \beta  = 0,$ $\beta \in \left( 0, \pi \right)$ with
summable potential $q \in L_{\mathbb{R}}^1 \left[0, \pi \right]$ is proved. This result is used to obtain
more precise asymptotic formulae for eigenvalues and norming constants of this problem.

\hfill Поступила 10 сентября 2017
\end{document}